\documentclass{article}

\usepackage{graphicx}                %
\title{Cobordism Theory and Poincar\'{e} Conjecture }      %
\author{Ming Yang(\tt{yangm926@163.com})\\  %
        Nanjing Agriculture University\\     %
        Nanjing, P.R.China 210032      %
        }   %
\date{}                     %
\begin{document}
\maketitle

\begin{abstract}                        %
In this paper, by use of techniques associated to cobordism theory
and Morse theory, we give a simple proof of Poincar\'{e} conjecture,
i.e.
Every compact smooth simply connected 3-manifold is homeomorphic to 3-sphere.                             %



\end{abstract}

\section{introduction}

Poincare conjecture, proposed by Henri Poincar\'{e} in 1904, is
presented in modern terminology as that every simply connected
closed 3-manifold without boundary is homeomorphic to the 3-sphere.

In 1961, by use of h-cobordism theory, S. Smale proved generalized
Poincar\'{e} conjecture: If $M^n$  is a differentiable homotopy
sphere of dimension $n\geq 5$ , then $M^n$  is homeomorphic to
$S^n$. However, it is huge difficulty to resolve Poincar\'{e}
conjecture of 3- and 4-dimension by using this theory. The reason is
that some results of h-cobordism theory cannot come into existence
in low dimensions, and all attempts to resolve these problems are
unsuccessful. They cannot prove that the right-hand spheres of
critical points with index 1 and the left-hand spheres of critical
points with index 2 of every simply connected, smooth, closed
3-manifold satisfy Morse first cancellation theorem and second
cancellation theorem. In this paper, by using some especial
techniques, we overcome these difficulties and give a proof of
Poincare conjecture. We prove that the second cancellation theorem
still hold in 3-manifolds, that is,  if the smooth compact
3-manifold $(W\cup W^\prime; S^2, S^2)$  is simply connected, $(W;
S^2, V)$ has one critical point of index 1 and $(W^\prime; V, S^2)$
has one critical point of index 2, then the intersection number
$S_L^1\cdot S_R^1=\pm1$ is true, there are at most finite smooth
isotopies such that $S_L^1$ and $S_R^1$ intersect at one point in
the oriented closed 2-manifold $V$, $W\cup W^\prime$ is
diffeomorphic to $S^2\times [0, 1]$. Moreover, we also prove the
claim that the right-hand spheres of critical points with index 1
and the left-hand spheres of critical points with index 2 of every
simply connected, smooth, closed 3-manifold satisfy the conditions
of Morse first cancellation theorem under at most finite smooth
isotopies.

 We shall now outline the contents of the paper. In section
2, we give the relationship between two linear expressions of
generators of $\pi_1(V)$, where $V$ is an oriented compact
2-manifold. In the linear expression, every entry of matrices is
determined uniquely by the intersection numbers among the
generators. In section 3, we prove that if $(W\cup W^\prime; S^2,
S^2)$ is simply connected then the coefficient matrix is a minor
diagonal matrix and every coefficient is 0 or $\pm1$. In section 4,
According to the arrangement theorem of elementary cobordisms and
the isotopy lemma, by use of proper arrangement of numbering of the
critical points, $S_L^1(q_i)$ and $S_R^1(p_i)$ intersect at one
point. In section 5, we obtain the conclusion that $(W\cup W^\prime;
S^2, S^2)$ is the product manifold by Morse's first cancellation
theorem.

\section{Linear Expression between Two Sets of Generators of
$\pi_1(V)$}

Let $T(k)$ be a smooth oriented closed 2-manifold with the genus $k$
and $g:S^1\rightarrow T(k)$ be a continuous mapping. $[g(S^1)]$
denotes the homotopy class of $g(S^1)$. If
$\alpha\cap\beta\neq\emptyset, \forall\alpha\in[g(S^1)],
\forall\beta\in[l(S^1)]$, there are two closed paths
$\alpha_0\in[g(S^1)]$ and $\beta_0\in[l(S^1)]$ with finite cross
points $x_1, x_2, \cdots, x_n$, each point is on both a smooth curve
segment of $\alpha_0$ and a smooth curve segment of $\beta_0$.

Suppose that $\omega$ is an orientation of $T(k)$. Given $S^1$ an
orientation, then, the orientations are fixed for $\alpha_0(S^1)$
and $\beta_0(S^1)$. $\alpha_0(S^1)$ has one tangent vector
$T(\alpha_0)_{x_i}$ at $x_i$, $\beta_0(S^1)$ has one tangent vector
$T(\beta_0)_{x_i}$ at $x_i$. When the orientation of the tangent
vector frame $(T(\alpha_0)_{x_i}, T(\beta_0)_{x_i})$ is the same as
the orientation $\omega$, the intersection number of $\alpha_0(S^1)$
and $\beta_0(S^1)$ at the point $x_i$ is 1, namely
$(\alpha_0\cdot\beta_0)_{x_i}=1$. When the orientation of the
tangent vector frame $(T(\alpha_0)_{x_i}, T(\beta_0)_{x_i})$ is the
same as the orientation $-\omega$, the intersection number of
$\alpha_0(S^1)$ and $\beta_0(S^1)$ at the point $x_i$ is $-1$,
namely $(\alpha_0\cdot\beta_0)_{x_i}=-1$.The intersection number of
$g(S^1))$ and $l(S^1))$ is defined as

$$l\cdot g=\sum{_x(l\cdot g)_x}$$

It is well known the intersection number $l\cdot g$ is a homotopy
invariant.

Let $l$  and $g$  be two oriented closed paths with a common point
$y$. $l\circ g$ denotes a oriented closed path starting at $y$
running back to $y$  along $g$ and again running back to $y$ along
$l$. $l^{-1}$ denotes the reverse of closed path $l$.

{\bf Definition 1.}\ Let $l$ and $g$  be two closed paths in $T(k)$.
There are $\alpha_0\in[g(S^1)]$ and $\beta_0\in[l(S^1)]$ which have
the least cross points, the number of the cross points of $\alpha_0$
and $\beta_0$ is called intersection degree, denoted by $d(l, g)$.

The intersection numbers and the intersection degrees have following
properties

(1)\ $d(l, g)=d(g, l)$ is a homotopy invariant;

(2)\ The intersection degree $d(l, g)$ is a nonnegative integer;

(3)\ There are $\alpha_0\in[g(S^1)]$ and $\beta_0\in[l(S^1)]$
satisfying $\alpha_0\cap\beta_0=\emptyset$, if and only if $d(l,
g)=0$;

(4)\ $l\cdot g=-g\cdot l,\ l^{-1}\cdot g=-l\cdot g$;

(5)\ $(l_1\circ l_2)\cdot g=l_1\cdot g+l_2\cdot g,\ g\cdot(l_1\circ
l_2)=g\cdot l_1+g\cdot l_2$;

(6)\ $l\cdot l=0$;

(7)\ $d(l,g)\geq \mid l\cdot g\mid$;

(8)\ $d(l, g)=0\Rightarrow l\cdot g=0$,\ $l\cdot g=0\Rightarrow d(l,
g)=0$ or $d(l, g)$ is positive even number.

According to the well-known theory of the oriented differentiable
closed 2-manifold ([2]), on the oriented differentiable closed
2-manifold $T(k)$  with the genus $k$, there exist $2k$
1-sub-manifolds $\{\alpha_1, \beta_1, \cdots, \alpha_k, \beta_k\}$
which satisfies the followings

(1)\ $\alpha_i$ and $\beta_i$ transversally intersect at one point
$(i=1, \cdots, k)$;

(2)\ $\alpha_i\cap \beta_j=\emptyset;\
\alpha_i\cap\alpha_j=\emptyset;\ \beta_i\cap\beta_j=\emptyset\
(\forall i\neq j)$;

(3)\ If $k=1$, $\pi_1(T(1))$ generated by $\{[\alpha], [\beta]\}$ is
a commutative group;

(4)\ If $k\geq 2$, $\pi_1(T(k))$ is a non-commutative group which is
the quotient of the free group on the generators  $\{\alpha_1,
\beta_1, \cdots, \alpha_k, \beta_k\}$ modulo the normal subgroup
generated by the element

$$\prod_{i=1}^k[\alpha_i,\beta_i]$$
where $[\alpha_i,\beta_i]=\alpha_i\beta_i\alpha_i^{-1}\beta_i^{-1}$.

Given $T(k)$ an orientation, and given $\alpha_i,\ \beta_i\ (i=1,
\cdots, k)$ proper orientations, we have

\begin{equation}
\alpha_i\cdot\beta_i=1\ (i=1,\cdots,k),\ \alpha_i\cdot\beta_j=0\
(i\neq j)
\end{equation}

\begin{equation}
\alpha_i\cdot\alpha_j=\beta_i\cdot\beta_j=0\ (\forall i, j)
\end{equation}

For $k\geq2$, $\forall \delta, \varepsilon\in \pi_1(T(k)),\ [\delta,
\varepsilon]=\delta\varepsilon\delta^{-1}\varepsilon^{-1}$ is called
the commutant of $\delta$ and $\varepsilon$. The normal subgroup
$[\pi_1(T(k)), \pi_1(T(k))]$ generated by all commutators is called
commutant subgroup.

{\bf Lemma 1.}\ if $\forall e\in [\pi_1(T(k)),\pi_1(T(k))]$ and
$\forall l\in \pi_1(T(k))$, then $e\cdot l=0$.

{\bf Proof.}\ $\forall \delta, \varepsilon\in \pi_1(T(k))$, $l\in
\pi_1(T(k))$, we have

$[\delta,\varepsilon]\cdot
l=(\delta\varepsilon\delta^{-1}\varepsilon^{-1})\cdot l=\delta\cdot
l+\varepsilon\cdot l+\delta^{-1}\cdot l+\varepsilon^{-1}\cdot
l=\delta\cdot l+\varepsilon\cdot l-\delta\cdot l-\varepsilon\cdot
l=0$

So, the conclusion is obtained by the definition of
$[\pi_1(T(k)),\pi_1(T(k))]$.\ {\bf QED}

Since the quotient group

$$\pi_1(T(k))/[\pi_1(T(k)), \pi_1(T(k))]$$
is a commutative group and
$\{\alpha_1,\beta_1,\alpha_2,\beta_2,\cdots,\alpha_k,\beta_k\}$ are
the infinite order generators, each element $l\in \pi_1(T(k))$  can
be expressed linearly as following

\begin{equation}
l=\sum_{i=1}^k m_i\alpha_i+\sum_{i=1}^k n_i\beta_i\pmod{
[\pi_1(T(k)),\pi_1(T(k))]}
\end{equation}
where, $m_i, n_i\ (1, \cdots, k)$ are all integers.

According to (1), (2) and Lemma 1, we have follows

$$m_i=l\cdot \beta_i,\ \ n_i=-l\cdot \alpha_i$$

(3) is written by follows

\begin{equation}
l=\sum_{i=1}^k (l\cdot \beta_i)\alpha_i-\sum_{i=1}^k (l\cdot
\alpha_i)\beta_i,\pmod{[\pi_1(T(k)),\pi_1(T(k))]}
\end{equation}

{\bf Lemma 2.}\ For $\forall l\in \pi_1(T(k))$, the linear
representation (4) is unique.

{\bf Proof.}\ As the numbers $l\cdot \alpha_i,\ l\cdot \beta_i\
(i=1,2,\cdots,k)$ are all the homotopy invariants, so every
coefficient in (4) is uniquely determined by $l$ and the generators
${\alpha_1,\beta_1,\alpha_2,\beta_2\cdots,\cdots,\alpha_k,\beta_k}$,
the linear representation (4) of $l$ is unique. {\bf QED}

Taking an element $g\in \pi_1(T(k))$. $g$ can be expressed linearly
as follows

$$g=\sum_{i=1}^k (g\cdot\beta_i)\alpha_i-\sum_{i=1}^k (g\cdot
\alpha_i)\beta_i,\pmod{[\pi_1(T(k)), \pi_1(T(k))]}$$

The intersection number $l\cdot g$ of $l$ and $g$ is

$$\begin{array}{l}
 l \cdot g = (\sum\limits_{i = 1}^k {(l \cdot \beta _i )\alpha _i } -
\sum\limits_{i = 1}^k {(l \cdot \alpha _i )\beta _i } ) \cdot
(\sum\limits_{i = 1}^k {(g \cdot \beta _i )\alpha _i } -
\sum\limits_{i =
1}^k {(g \cdot \alpha _i )\beta _i } ) \\
 = - \sum\limits_{i = 1}^k {(l \cdot \beta _i )} (g \cdot \alpha _i ) +
\sum\limits_{i = 1}^k {(l \cdot \alpha _i )} (g \cdot \beta _i ) \\
 = \sum\limits_{i = 1}^k {((l \cdot \alpha _i )(g \cdot \beta _i ) - (l
\cdot \beta _i )} (g \cdot \alpha _i )) \\
 = \sum\limits_{i = 1}^k {\det \left( {{\begin{array}{*{20}c}
 {l \cdot \beta _i } \hfill & { - l \cdot \alpha _i } \hfill \\
 {g \cdot \beta _i } \hfill & { - g \cdot \alpha _i } \hfill \\
\end{array} }} \right)} \\
 \end{array}$$

Let $\{\alpha_1,\beta_1,\cdots,\alpha_k,\beta_k\}$ be the generators
of $\pi_1(T(k))$ and $h:T(k)\rightarrow T(k)$ be a diffeomorphism.
If $h$ gives $T(k)$ a opposite orientation $h^\ast(\omega)=-\omega$,
$h$ is called the cobordism diffeomorphism.

If $h:T(k)\rightarrow T(k)$ is a cobordism diffeomorphism, then
$\{h(\alpha_1), h(\beta_1), \cdots, h(\alpha_k), h(\beta_k)\}$ are
all 1-submanifolds which have the following properties

(1)\ $h(\alpha_i)$ and $h(\beta_i)$\ $(i=1, \cdots, k)$
transversally intersect at one point;

(2)\ $h(\alpha_i)\cap h(\beta_j)=h(\alpha_i)\cap
h(\alpha_j)=h(\beta_i)\cap h(\beta_j)=\emptyset$\ $(\forall i\neq
j)$;

(3)\ If $k=1$, $\pi_1(T(1))$ generated by $\{[h(\alpha)],
[h(\beta)]\}$ is a commutative group;

(4)\ If $k\geq 2$, $\pi_1(T(k))$ is a non-commutative group. It is
the quotient of the free group on the generators  $\{h(\alpha_1),
h(\beta_1), \cdots, h(\alpha_k), h(\beta_k)\}$ modulo the normal
subgroup generated by the element

$$\prod_{i=1}^k[h(\alpha_i), h(\beta_i)]$$

Given $T(k)$ the orientation $-\omega$ and given $\{h(\alpha_i),\
h(\beta_i) \mid i=1,\cdots,k\}$ the proper orientations, the
equations can be obtained

$$h(\alpha_i)\cdot h(\beta_i)=1\ (i=1,\cdots,k);\ h(\alpha_i)\cdot h(\beta_j)=0\
(\forall i\neq j)$$
$$h(\alpha_i)\cdot h(\alpha_j)=h(\beta_i)\cdot h(\beta_j)=0\ (\forall i, j=1,\cdots,k)$$

Now use the marks $\theta_i=h(\alpha_i),\ \gamma_i=h(\beta_i)\
(i=1,\cdots,k)$ and

$$(\alpha,
\beta)=\{\alpha_1,\cdots,\alpha_k, \beta_1,\cdots,\beta_k\};\
(\theta, \gamma)=\{\theta_1,\cdots,\theta_k,
\gamma_1,\cdots,\gamma_k\}$$

By use of $(\alpha, \beta)$ to express linearly $(\theta, \gamma)$,
we have follows.

\begin{equation}
\theta_i=\sum_{j=1}^k (\theta_i\cdot \beta_j)\alpha_j-\sum_{j=1}^k
(\theta_i\cdot \alpha_j)\beta_j,\pmod{[\pi_1(T(k)),\pi_1(T(k))]}
\end{equation}

\begin{equation}
 \gamma_i=\sum_{j=1}^k (\gamma_i\cdot
\beta_j)\alpha_j-\sum_{j=1}^k (\gamma_i\cdot
\alpha_j)\beta_j,\pmod{[\pi_1(T(k)),\pi_1(T(k))]}
\end{equation}

On the other hand, using $(\theta,\gamma)$ to express linearly
$(\alpha,\beta)$,

\begin{equation}
\alpha_i=\sum_{j=1}^k (\alpha_i\cdot \gamma_j)\theta_j-\sum_{j=1}^k
(\alpha_i\cdot \theta_j)\gamma_j,\pmod{[\pi_1(T(k)),\pi_1(T(k))]}
\end{equation}

\begin{equation}
\beta_i=\sum_{j=1}^k (\beta_i\cdot \gamma_j)\theta_j-\sum_{j=1}^k
(\beta_i\cdot \theta_j)\gamma_j,\pmod{[\pi_1(T(k)),\pi_1(T(k))]}
\end{equation}

In order to we rewrite (5), (6), (7), (8) in terms of matrices.
Define

$$\alpha=(\alpha_1,\cdots,\alpha_k);\ \ \beta=(\beta_1,\cdots,\beta_k);$$
$$\theta=(\theta_1,\cdots,\theta_k);\ \ \gamma=(\gamma_1,\cdots,\gamma_k);$$

$$
\theta ^T \cdot \alpha = \left( {{\begin{array}{*{20}c}
 {\theta _1 \cdot \alpha _1 } \hfill & \cdots \hfill & {\theta _1 \cdot
\alpha _k } \hfill \\
 \vdots \hfill & \cdots \hfill & \vdots \hfill \\
 {\theta _k \cdot \alpha _1 } \hfill & \cdots \hfill & {\theta _k \cdot
\alpha _k } \hfill \\
\end{array} }} \right) = \left( {\theta _i \cdot \alpha _j }
\right)_{k\times k}
$$

and\ $\theta^T\cdot\beta=(\theta_i\cdot\beta_j)_{k\times k}$;\ \
$\gamma^T\cdot\alpha=(\gamma_i\cdot\alpha_j)_{k\times k}$;\ \
$\gamma^T\cdot\beta=(\gamma_i\cdot\beta_j)_{k\times k}$.

(5) and (6) can be written as

\begin{equation}
\left( {{\begin{array}{*{20}c}
 {\theta ^T} \hfill \\
 {\gamma ^T} \hfill \\
\end{array} }} \right) = \left( {{\begin{array}{*{20}c}
 {\theta ^T \cdot \beta } \hfill & { - \theta ^T \cdot \alpha } \hfill \\
 {\gamma ^T \cdot \beta } \hfill & { - \gamma ^T \cdot \alpha } \hfill \\
\end{array} }} \right)\left( {{\begin{array}{*{20}c}
 {\alpha ^T} \hfill \\
 {\beta ^T} \hfill \\
\end{array} }} \right) = H\left( {{\begin{array}{*{20}c}
 {\alpha ^T} \hfill \\
 {\beta ^T} \hfill \\
\end{array} }} \right)\quad\pmod{[\pi _1 (T(k)),\,\pi _1 (T(k))]}
\end{equation}

(7) and (8) can be written as

\begin{equation}
\left( {{\begin{array}{*{20}c}
 {\alpha ^T} \hfill \\
 {\beta ^T} \hfill \\
\end{array} }} \right) = \left( {{\begin{array}{*{20}c}
 {\alpha ^T \cdot \gamma } \hfill & { - \alpha ^T \cdot \theta } \hfill \\
 {\beta ^T \cdot \gamma } \hfill & { - \beta ^T \cdot \theta } \hfill \\
\end{array} }} \right)\left( {{\begin{array}{*{20}c}
 {\theta ^T} \hfill \\
 {\gamma ^T} \hfill \\
\end{array} }} \right) = H^{ - 1}\left( {{\begin{array}{*{20}c}
 {\theta ^T} \hfill \\
 {\gamma ^T} \hfill \\
\end{array} }} \right)\quad\pmod{[\pi _1 (T(k)),\,\pi _1 (T(k))]}
\end{equation}

As $\{\alpha_1,\cdots,\alpha_k, \beta_1,\cdots,\beta_k\}$ are the
generators of $\pi_1(T(k))$ and $\{\theta_1,\cdots,\theta_k,
\gamma_1,\cdots,\gamma_k\}$ also are the generators of
$\pi_1(T(k))$, therefore, $H$ and $H^{-1}$ are both the nonsingular
matrixes.

According the equations (9) and (10), we have

$$ \left( {{\begin{array}{*{20}c}
 {\theta ^T \cdot \beta } \hfill & { - \theta ^T \cdot \alpha } \hfill \\
 {\gamma ^T \cdot \beta } \hfill & { - \gamma ^T \cdot \alpha } \hfill \\
\end{array} }} \right)\left( {{\begin{array}{*{20}c}
 {\alpha ^T \cdot \gamma } \hfill & { - \alpha ^T \cdot \theta } \hfill \\
 {\beta ^T \cdot \gamma } \hfill & { - \beta ^T \cdot \theta } \hfill \\
\end{array} }} \right) = HH^{ - 1} = E_{2k}
$$ where $E_{2k}$ is a unit matrix.

Moreover, $H$ and $H^{-1}$ are both the integral matrixes, so

\begin{equation}
\det \left( {{\begin{array}{*{20}c}
 {\theta ^T \cdot \beta } \hfill & { - \theta ^T \cdot \alpha } \hfill \\
 {\gamma ^T \cdot \beta } \hfill & { - \gamma ^T \cdot \alpha } \hfill \\
\end{array} }} \right) = \det \left( {{\begin{array}{*{20}c}
 {\alpha ^T \cdot \gamma } \hfill & { - \alpha ^T \cdot \theta } \hfill \\
 {\beta ^T \cdot \gamma } \hfill & { - \beta ^T \cdot \theta } \hfill \\
\end{array} }} \right) = \pm 1
\end{equation}

For the determination of the coefficient matrixes $H$ and $H^{-1}$
in (9) and (10), it is necessary to discuss the cobordisms of the
oriented compact 3-manifolds.

\section{Linear Expression of Cobordism of Simply Connected
3-Manifolds}

{\bf Definition 2.}\ $(W;V_0,V_1)$ possessing a Morse function is
the triad of a smooth 3-manifold, if $W$ is a compact smooth
3-manifold with the boundaries $BdW$ which are disjoint union of two
both open and closed 2-manifolds $V_0$ and $V_1$.

$(W;V_0,V_1)$ is the elementary cobordism, if $W$ has exactly one
critical point.

$S^{n-1}$ represents the boundary of the unit closed disk $D^n$ in
$R^n$; and $OD^n$ represents the unit open disk in $R^n$.

Let $(W;V_0,V_1)$ be an elementary cobordism with Morse function
$f:W\rightarrow R^1$ and gradient vector field $\xi$ for $f$.
Suppose $p\in W$ is a critical point, and $V_0=f^{-1}(0)$ and
$V_1=f^{-1}(1)$ such that $0<f(p)<1$ and $c=f(p)$ is the only
critical value in the interval $[0,1]$. Define the left-hand sphere
$S_L$ of $p$ is just the intersection of $V_0$ with all integral
curves of $\xi$ leading to the critical point $p$. The left-hand
disc $D_L$ is a smoothly imbedded disc with boundary $S_L$, defined
to be the union of the segments of these integral curves beginning
in $S_L$ ending at $p$. The right-hand sphere $S_R$ of $p$ in $V_1$
is the boundary of segments of integral curves of $\xi$ beginning at
$p$ and ending in $S_R$.

{\bf Theorem 1.}\ Let $(W;V_0,V)$ be a triad of an oriented smooth
compact 3-manifold, and $V_0$ be diffeomorphic to $S^2$.
$f:W\rightarrow R^1$ is a Morse function, $f^{-1}(0)=V_0,
f^{-1}(1)=V$. There exist $k$ critical points $p_1,p_2,\cdots,p_k;
\lambda(p)=1$ in $W$, they are on one same horizontal plane. Then,

(1)\ $V$ is an oriented compact 2-manifold with the genus $k$.

(2)\ In $V$,  there exist $k$ 1-submanifolds
$\{\beta_1,\cdots,\beta_k\}$, which are not mutually intercrossed.

The right-hand spheres $S_R^1(p_1), \cdots, S_R^1(p_k)$ of the
critical points are also non-crossing 1-submanifolds in $V$.
$\beta_i$ only intercrosses with $S_R^1(p_i)$ at one point. The
homotopy classes $\{[S_R^1(p_i)],[\beta_i] | i=1,2,\cdots,k\}$ are
generators of $\pi_1(T(k))$,

(3)\ $\{[\beta_i] | i=1,2,\cdots,k\}$ is a set of generators of
$\pi_1(W)$, and $\pi_1(W)$ is the free product of $k$ infinite
cyclic groups.

{\bf Proof.}\ $W$ is an oriented smooth 3-manifold, so $V$ is an
oriented smooth 2-manifold. In $W$, there are $k$ non-crossing
characteristic embeddings ([1] P28)

$$\varphi_i: S^0\times OD^2\rightarrow V_0, (i=1,2,\cdots,k)$$

Taking a the disjoint sum

$$(V_0-\sum_{i=1}^k \varphi_i(S^0\times0))+(OD^1\times
S^1)_1+\cdots+(OD^1\times S^1)_k$$ and the equivalents as follows:

$$\varphi_i(u,\theta v)\sim (\theta u, v)_i,\ \ u\in S^0,\ v\in S^1,\ 0<\theta<1$$

Thus we obtain an oriented smooth 2-manifold $\chi(V_0,
\varphi_1,\cdots, \varphi_k)$ with genus is $k$, and
$S_R^1(p_i)=(0\times S^1)_i$, $S_L^0(p_i)=\varphi_i(S^0\times 0)$.

According to Theorem 3.13 of ([1]P31-36), $(W;V_0,V)$ is
diffeomorphic to
$$(\omega(V_0,\varphi_1,\cdots,\varphi_k), V_0,
\chi(V_0, \varphi_1,\cdots, \varphi_k))$$

Therefore, $V=\chi(V_0, \varphi_1,\cdots, \varphi_k))$.

Take $k$ 1-submanifolds $\{\beta_1,\cdots, \beta_k\}$ in $\chi(V_0,
\varphi_1,\cdots, \varphi_k)$ as follows:

First of all, we take non-crossing closed 2-disk $B_1,\cdots,B_k$ in
$V_0$, such that $\varphi_i(S^0\times OD^2)\subset B_i$. The genus
of $\chi(B_i,\varphi_i)$ is equal to 1,
$Bd\chi(B_i,\varphi_i)=BdB_i$. So, there exist a diffeomorphim
$\chi(B_i,\varphi_i)\cong S^1\times S^1-B_0$, here $B_0$ is an open
2-disk, in this diffeomorphism, $S_R^1(p_i)=(0\times S^1)_i$
corresponds to $y\times S^1\subset S^1\times S^1-B_0$. We retake one
1-submanifold $S^1\times z$ in $(S^1\times S^1-B_0)$ such that
$y\times S^1$ and $S^1\times z$ transversely intersect at one point
$y\times z$. In $\chi(B_i,\varphi_i)$, the 1-submanifold
corresponding to $S^1\times z$ is denoted by $\beta_i$, then
$S_R^1(p_i)$ and $\beta_i$ intersect just at one point, and two
homotopy classes $\{[\beta_i], [S_R(p_i)]\}$ generate
$\pi_1(\chi(B_i,\varphi_i))$.

Let $H_k=V_0-\cup_{i=1}^k(IntB_i)$, we have

$$\chi(V_0,\varphi_1,\cdots,\varphi_k)=H_k\cup\chi(B_1,\varphi_1)\cup\cdots\cup\chi(B_k,\varphi_k)$$

$$S_R(p_i)\cup\beta_i\subset Int\chi(B_i,\varphi_i),
(i=1,\cdots,k)$$

According to well-known 2-manifolds theory, the homotopy classes
$\{[\beta_i], [S_R(p_i)]\mid i=1,\cdots, k\}$ are the generators of
$\pi_1(V)$.

Let $\Psi:V_0\rightarrow BdD^3$ be a diffeomorphism, then
$M=W\cup_\Psi D^3$ is a smooth oriented 3-manifold. $M$  has a
deformation retract ([1]Theorem 3.14)

$$D^3\cup D_L^1(p_1)\cup\cdots\cup D_L^1(p_k)$$ where $D_L^1(p_i)$
are disjoint 1-discs, $D^3\cap
D_L^1(p_i)=S_L^0(p_i)=\varphi_i(S^0\times0)$.

As $D^3$ has a deformation retract to the origin, $M=W\cup_\Psi D^3$
has a deformation retract which is $k$ circles with one common
point. Hence, $\pi_1(W)=\pi_1(M)$ is the free products of $k$
infinite cyclic groups.

Since $S_R^1(p_i)$ is the boundary of 2-disc $D_R^2(p_i)$ and
$D_R^2(p_i)\subset W$, $S_R^1(p_i)$ is null homotopy in $W$.
According to the definition of $\{\beta_i \mid i=1,\cdots,k\}$,
$\{[\beta_i] | i=1,\cdots,k\}$ is a set of generators of
$\pi_1(W)$.\ {\bf QED}

{\bf Theorem 2.}\ Suppose that $(W; S^2,V)$ and $(W'; V, S^2)$ be
two oriented smooth compact 3-manifolds. $(W; S^2,V)$ has exactly
the critical points $p_1,\cdots,p_k\ (k\geq1)$ of type 1 and
$(W^\prime; V,S^2)$ has exactly the critical points $q_1,\cdots,q_k$
of type 2, then there exists a diffeomorphism $G:(W;
S^2,V)\rightarrow (W^\prime; S^2,V)$ such that

$$G(D_L^1(p_i))=D_R^1(q_i);\  G(D_R^2(p_i))=D_L^2(q_i)\
(i=1,\cdots,k)$$
$$G(V)=V;\ G(S^2)=S^2$$

{\bf Proof.}\ Let $\varphi_i: S^0\times OD^2\rightarrow S^2$ be
disjoint embeddings of $p_i$ and $\rho_i: S^1\times OD^1\rightarrow
V$ be disjoint embeddings of $q_i$.

Since each $q_i$ is the critical point of index 1 in $(W'; S^2,V)$,
$(\rho_i)_R: S^0\times OD^2\rightarrow S^2$ are disjoint embeddings
of $q_i$ in $(W^\prime; S^2,V)$. Then there exists an isotopy $h:
S^2\times I\rightarrow S^2$ such that

$$h_1(\varphi_i(S^0\times 0))=(\rho_i)_R(S^0\times 0)\
(i=1,\cdots, k)$$

It is known ([1] Theorem 3.13) that $(W; S^2, V)$ is diffeomorphic
to
$$(\omega(S^2,\ \varphi_1,\cdots,\ \varphi_k);\ S^2,\ \chi(S^2,\ \varphi_1,
\cdots,\ \varphi_k
))$$ and $(W^\prime;\ S^2,\ V)$ is diffeomorphic to

$$(\omega(S^2,\ (\rho_1)_R,\cdots,\ (\rho_k)_R);\
S^2,\ \chi(S^2,\ (\rho_1)_R,\cdots,\ (\rho_k)_R))$$

Hence $h_1$ can be extended to the diffeomorphism $G:(W;\ S^2,\
V)\rightarrow (W^\prime;\ S^2,\ V)$ which satisfies following

$$G(D_L^1(p_i))=D_R^1(q_i);\ \ G(D_R^2(p_i))=D_L^2(q_i)\ (i=1,\cdots,k)$$
$$G(V)=V;\ G(S^2)=S^2$$
\ \ {\bf QED}

Suppose that $(W; V_1,V)$, $(W^\prime; V,V_2)$, $(W\cup W^\prime;
V_1,V_2)$ satisfy the conditions in Theorem 2.

Given $W$ an orientation $\xi$, $V\subset BdW$ has a induced
orientation $\omega$. The frame of the tangent vector
$(\tau_1,\tau_2)$ at some point $x\in V$ of $V$ is positively
oriented, if the 3-frame $(\nu,\tau_1,\tau_2)$ is positively
oriented in $TW_x$, where $\nu$ is any vector at $x$ tangent to $W$
but not to $V$ and pointing out of $W$.

The diffeomorphism $G$ gives an orientation  $G^\ast(\xi)$ of
$W^\prime$. There exists the unique diffeomorphism $h: V\rightarrow
V$ such that
$$W\cup_{h\circ G_V} W^\prime=W\cup W^\prime$$
$$h(S_R^1(p_i))=S_L^1(q_i),(i=1,\cdots,k)$$

So $h: V\rightarrow V$ gives $V$ an opposite orientation
$$h^\ast(\omega)=-\omega$$

Since $\{[S_R^1(p_i)],[\beta_i]\mid i=1,2,\cdots,k\}$ is the set of
generators of $\pi_1(V)$,
$$\{[h(S_R^1(p_i))],[h(\beta_i)]\mid
i=1,2,\cdots,k\}=\{[S_L^1(q_i )],[h(\beta_i)]\mid i=1,\cdots,k\}$$
is also the set of generators of $\pi_1(V)$. Hence $h:V\rightarrow
V$ is a cobordism diffeomorphim.

{\bf Definition 3.}\ Let
$\{\alpha_1,\beta_1,\cdots,\alpha_k,\beta_k\}$ be the generators of
$\pi_1(V)$ and $V$ be an oriented compact 2-manifold. Each element
$g\in \pi_1(V)$ is expressed in the product form

\begin{equation}
g=\xi_1\cdots\xi_m,\ \
\xi\in\{\alpha_1^{\pm1},\beta_1^{\pm1},\cdots,\alpha_k^{\pm1},\beta_k^{\pm1}\}
\end{equation}

$\alpha_j^{\pm1}$\ (or $\beta_k^{\pm1}$) is called as homogenous, if
the numbers of occurrence $\alpha_j$\ ($\beta_j$) and
$\alpha_j^{-1}$\ ($\beta_j^{-1}$) are the same in equation (12).

{\bf Lemma 3.}\ Suppose that $g\in \pi_1(V)$ and $g$ be expressed in
the product form (12). $\alpha_j^{\pm1}$ are homogenous if and only
if $g\cdot\beta_j=0$. $\beta_j^{\pm1}$ are homogenous if and only if
$g\cdot\alpha_j=0$.

{\bf Proof.}\ Suppose that $\alpha_j^{\pm1}$ are homogenous in (12).
We have follows.
$$\alpha_i\cdot\beta_i=+1;\  \alpha_i\cdot\beta_j=0\ \ (i\neq j)$$
$$\alpha_i\cdot\alpha_j=\beta_i\cdot\beta_j=0\ \ (\forall
i,j=1,\cdots,k)$$

So $g\cdot\beta_j=0$. If the numbers of occurrence $\alpha_j$ and
$\alpha_j^{-1}$ are not the same in Equation (12), then it is
obvious that $g\cdot\beta_j\neq 0$, hence, if $g\cdot\beta_j=0$,
$\alpha_j^{\pm1}$ are homogenous in (12).

The second conclusion will be obtained by the same reason.{\bf QED}

{\bf Lemma 4.}\ Suppose that $(W; V_1, V)$, $(W^\prime; V, V_2)$,
$(W\cup W^\prime; V_1, V_2)$ satisfy the conditions in Theorem 2.
Let $l\in\pi_1(V)$, if at least one of these $k$ integers
$l\cdot\alpha_1,\ l\cdot\alpha_2,\cdots,\ l\cdot\alpha_k$ is not
equal to 0, then, $l$ is not null homotopy in $W$. if at least one
of these $k$ integers $l\cdot\theta_1,\ l\cdot\theta_2,\cdots,\
l\cdot\theta_k$ is not equal to 0, then, $l$ is not null homotopy in
$W^\prime$.

{\bf Proof.}\ On the boundary $V$ of $W$,  by use of $(\alpha,
\beta)$ to represent $l$ linearly

\begin{equation}
l=\sum_{i=1}^k (l\cdot\beta_i)\alpha_i-\sum_{i=1}^k
(l\cdot\alpha_i)\beta_i,\pmod{[\pi_1(V), \pi_1(V)]}
\end{equation}
$l$ is expressed in the product form of

\begin{equation}
l=\xi_1\cdots\xi_n,\ \ \xi_i\in
\{\alpha_1^{\pm1},\beta_1^{\pm1},\cdots,\alpha_k^{\pm1},\beta_k^{\pm1}\}
\end{equation}

Because $\{\alpha_i \mid i=1,\cdots,k\}$ are all null homotopy in
$W$, so in $\pi_1(W)$, $l$ can be expressed as

\begin{equation}
l=\sigma_1\cdots\sigma_n,\ \ \sigma\in
\{\beta_1^{\pm1},\cdots,\beta_k^{\pm1}\}
\end{equation}

Where $\sigma_1\cdots\sigma_n$ is obtained by removing all
$\{\alpha_i,\alpha_i^{-1}\ \mid i=1,\cdots,k\}$ and keeping all
$\{\beta_i,\beta_i^{-1} \mid i=1,\cdots,k\}$ same order as that in
(14).

In (13), if at least one of these $k$ integers $l\cdot\alpha_1,\
l\cdot\alpha_2,\cdots,\ l\cdot\alpha_k$ is not equal to 0, for
example, $l\cdot\alpha_j\neq0$, then $\beta_j^{\pm1}$ are not
homogenous in the equation (14) and (15). Since $\pi_1(W)$ generated
by the homotopy classes $\{\beta_1,\cdots,\beta_k\}$ is the free
product of $k$ infinite cyclic groups, $l$ is not null homotopy in
$W$.

The second conclusion will be obtained from the diffeomorphisms
$G:(W^\prime; V_2,V) \rightarrow (W; V_1,V)$ and $h: V\rightarrow
V$.  {\bf QED}

{\bf Lemma 5.}\ Suppose that $(W; V_1,V)$, $(W^\prime; V,V_2)$,
$(W\cup W^\prime; V_1,V_2)$ are three oriented smooth compact
3-manifolds and that $W\cup W^\prime,\ V_1,\ V_2$ are all simply
connected. $f:(W\cup W^\prime; V_1,V_2)\rightarrow R^1$ is a Morse
function with the critical points $p_1,\cdots,p_k\ (k\geq1)$ of type
1 and the critical points $q_1,\cdots,q_k$ of type 2,
$f^{-1}(-2)=V_1$, $f^{-1}(0)=V$, $f^{-1}(2)=V_2$, $f(p_i)=-1$,
$f(q_i)=+1,\ (i=1,\cdots,k)$. then, the homotopy classes
$\{\alpha_1,\theta_1,\cdots,\alpha_k,\theta_k\}$ are the generators
of $\pi_1(V)$. And

$$\gamma^T\cdot\alpha=0;\ \beta^T\cdot\theta=0$$

{\bf Proof.}\ If $\gamma^T\cdot\alpha\neq0$, then, there exists a
non-zero row vector $(\gamma_j\cdot\alpha_1,\
\gamma_j\cdot\alpha_2,\cdots,\ \gamma_j\cdot\alpha_k)$ in the matrix
$\gamma^T\cdot\alpha$. Assuming that there is at least one non-zero
number in set $\{\gamma_j\cdot\alpha_1,\
\gamma_j\cdot\alpha_2,\cdots,\ \gamma_j\cdot\alpha_k\}$. Adopting
the linear representation of $\gamma_j$,

$$\gamma_j=\sum_{i=1}^k (\gamma_j\cdot\beta_i)\alpha_i-\sum_{i=1}^k
(\gamma_j\cdot\alpha_i)\beta_i,\pmod{[\pi_1(V), \pi_1(V)]}$$

and a product representation

$$\gamma_j=\sigma_1\cdots\sigma_n,\ \sigma\in\{\beta_1^{\pm1},\cdots,\beta_k^{\pm1}\}$$

According to Lemma 4, $\gamma_j$ is not null homotopy in $W$.

Taking a positive number $\varepsilon$, $\varepsilon<1$, and two
path connected sets $X=W\cup f^{-1}[0,\varepsilon)$, $Y=W^\prime\cup
f^{-1}(-\varepsilon,0]$, then, $X\cap Y=f^{-1}(-\varepsilon,
\varepsilon)$. $f^{-1}(-\varepsilon, \varepsilon)$ has not any
critical point, so $f^{-1}(-\varepsilon, \varepsilon)$ is an open
product manifold $V\times(-\varepsilon, \varepsilon)$. Since of $X$,
$Y$, $X\cap Y$ are all path connected open sets, $\{X,\ Y,\ X\cap
Y\}$ is a path connected open covering of $W\cup W^\prime$.

$W$ is a deformation retract of $\overline{X}=W\cup f^{-1}[0,
\varepsilon]$ and $Y$ is a deformation retract of
$\overline{Y}=W^\prime\cup f^{-1}[-\varepsilon,0]$; $\gamma_j$ is
not null homotopy in $W$ and $W^\prime$. Combining Lemma 4, Theorem
1 with Van. Kampen theorem, the conclusion $\pi_1(W\cup
W^\prime)\neq1$ is obtained. However, it is known that $\pi_1(W\cup
W^\prime)=1$, hence, $\gamma^T\cdot\alpha=0$.

$\beta^T\cdot\theta=0$ is obtained by use of the same technique.\
{\bf QED}

If $W\cup W^\prime$ satisfy the conditions in Lemma 5, according to
Lemma 5, we have $\gamma^T\cdot\alpha=0$ and $\beta^T\cdot\theta=0$,
so $\alpha^T\cdot\gamma=0$; $\theta^T\cdot\beta=0$. In (9), (10),
the square matrixes of coefficients are both nonsingular matrixes,
so $-\theta^T\cdot\alpha,\ \gamma^T\cdot\beta,\
-\alpha^T\cdot\theta,\ \beta^T\cdot\gamma$ are all nonsingular
matrixes. Hence we have

\begin{equation}
\left( {\begin{array}{l}
 \theta ^T \\
 \gamma ^T \\
 \end{array}} \right) = \left( {{\begin{array}{*{20}c}
 0 \hfill & { - \theta ^T \cdot \alpha } \hfill \\
 {\gamma ^T \cdot \beta } \hfill & 0 \hfill \\
\end{array} }} \right)\left( {\begin{array}{l}
 \alpha ^T \\
 \beta ^T \\
 \end{array}} \right)\quad\pmod{[\pi_1(V), \pi_1(V)]}
\end{equation}

\begin{equation}
\left( {\begin{array}{l}
 \alpha ^T \\
 \beta ^T \\
 \end{array}} \right) = \left( {{\begin{array}{*{20}c}
 0 \hfill & { - \alpha ^T \cdot \theta } \hfill \\
 {\beta ^T \cdot \gamma } \hfill & 0 \hfill \\
\end{array} }} \right)\left( {\begin{array}{l}
 \theta ^T \\
 \gamma ^T \\
 \end{array}} \right)\quad\pmod{[\pi_1(V), \pi_1(V)]}
\end{equation}

For convenient, we will always use the following expressions.
$$\alpha_i=S_R^1(p_i);\ \theta_i=S_L^1(q_i)=h(S_R^1(p_i));\ \gamma_i=h(\beta_i),\
(i=1,\cdots,k)$$

{\bf Theorem 3.}\ Suppose that $(W; V_1, V)$, $(W^\prime; V, V_2)$,
$(W\cup W^\prime; V_1, V_2)$ satisfy the conditions in Lemma 5.
Then, $\{\alpha_i, \beta_i, \theta_i, \gamma_i \mid i=1,\cdots,k\}$
are all 1-submanifolds in V and the homotopy classes $\{\alpha_1,
\theta_1,\cdots, \alpha_k, \theta_k\}$ are the generators of
$\pi_1(V)$. $\{\alpha_i, \beta_i \mid i=1,\cdots,k\}$ satisfy the
conditions in theorem 1. $\pi_1(W^\prime)$ is isomorphic to
$\pi_1(W)$, $\pi_1(W^\prime)$ generated by $\{\gamma_1,\cdots,
\gamma_k\}$ is the free product of $k$ infinite cyclic groups.
Moreover, we have

\begin{equation}
\theta_i=\pm\beta_{\sigma(i)}\pmod{[\pi_1(V), \pi_1(V)]}
\end{equation}

\begin{equation}
\gamma_i=\pm\alpha_{\sigma(i)}\pmod{[\pi_1(V), \pi_1(V)]}
\end{equation}
where $\sigma$ is a permutation of $\{1,\cdots,k\}$

{\bf Proof.}\ According to Theorem 1, there are 1-submanifolds
$\{\alpha_i, \beta_i \mid i=1,\cdots,k\}$ in $V$, they are the
generators of $\pi_1(V)$ and $\{\beta_i \mid i=1,\cdots,k\}$ are the
generators of $\pi_1(W)$ which is the free product of infinite cycle
groups. From the cobordism diffeomorphism $h: V\rightarrow V$ and
the diffeomorphism $G: (W; S^2,V)\rightarrow(W^\prime; S^2,V)$, we
obtain that $\pi_1(V)$ generated by $\{\theta_i, \gamma_i \mid
i=1,\cdots,k\}$ which are all 1-submanifolds in $V$ and
$\pi_1(W^\prime)$ generated by $\{\gamma_1,\cdots, \gamma_k\}$ is
the free product of $k$ infinite cyclic groups. Moreover, since
$\{\alpha_1,\cdots, \alpha_k\}$ are disjoint 1-submanifolds,
$\{\theta_1,\cdots, \theta_k\}$ are also disjoint 1-submanifolds.

We have

$$\alpha_i\cap\beta_i\ \mbox{is just one point};\ d(\alpha_i, \beta_i)=1;\ \alpha_i\cdot\beta_i=1$$
$$\theta_i\cap\gamma_i\ \mbox{is just one point};\ d(\theta_i\cdot\gamma_i)=1;\ \theta_i\cdot\gamma_i=-1$$
$$\alpha_i\cap\beta_j=\alpha_i\cap\alpha_j=\beta_i\cap\beta_j=\emptyset,\ (\forall i\neq j)$$
$$\theta_i\cap\gamma_j=\theta_i\cap\theta_j=\gamma_i\cap\gamma_j=\emptyset,\ (\forall i\neq j)$$
and

$$\begin{array}{l}
 \theta_l \cdot \gamma_j = (\sum\limits_{i = 1}^k {(\theta_l \cdot \beta _i )\alpha _i } -
\sum\limits_{i = 1}^k {(\theta_l \cdot \alpha _i )\beta _i } ) \cdot
(\sum\limits_{i = 1}^k {(\gamma_j \cdot \beta _i )\alpha _i } -
\sum\limits_{i =
1}^k {(\gamma_j \cdot \alpha _i )\beta _i } ) \\
 = - \sum\limits_{i = 1}^k {(\theta_l \cdot \beta _i )} (\gamma_j \cdot \alpha _i ) +
\sum\limits_{i = 1}^k {(\theta_l \cdot \alpha _i )} (\gamma_j \cdot \beta _i ) \\
 = \sum\limits_{i = 1}^k {((\theta_l \cdot \alpha _i )(\gamma_j \cdot \beta _i ) - (\theta_l
\cdot \beta _i )} (\gamma_j \cdot \alpha _i )) \\
 = \sum\limits_{i = 1}^k {\det \left( {{\begin{array}{*{20}c}
 {\theta_l \cdot \beta _i } \hfill & { - \theta_l \cdot \alpha _i } \hfill \\
 {\gamma_j \cdot \beta _i } \hfill & { - \gamma_j \cdot \alpha _i } \hfill \\
\end{array} }} \right)} \\
 \end{array}$$

\begin{equation}
\begin{array}{l} d(\theta_l, \gamma_j)\geq \sum\limits_{i=1}^k
{|\det \left( {{\begin{array}{*{20}c}
 {\theta_l \cdot \beta _i } \hfill & { - \theta_l \cdot \alpha _i } \hfill \\
 {\gamma_j \cdot \beta _i } \hfill & { - \gamma_j \cdot \alpha _i } \hfill \\
\end{array} }} \right)|} \\
 \end{array}
\end{equation}

From (16) and (17), we obtain

$$\theta^T=-(\theta^T\cdot\alpha)\beta^T \quad\pmod{[\pi_1(V),
\pi_1(V)]}$$

$$\gamma^T=(\gamma^T\cdot\beta)\alpha^T \quad\pmod{[\pi_1(V),
\pi_1(V)]}$$

According to Lemma 5 and (20), we obtain

$$\begin{array}{l}
1=d(\theta_j, \gamma_j)\geq \sum\limits_{i=1}^k {|\det
\left({{\begin{array}{*{20}c}
 {\theta_j \cdot \beta _j } \hfill & { - \theta_j \cdot \alpha _i } \hfill \\
 {\gamma_j \cdot \beta _i } \hfill & { - \gamma_j \cdot \alpha _i } \hfill \\
\end{array} }} \right)|}=\sum\limits_{i=1}^k {|\det
\left({{\begin{array}{*{20}c}
 {\ \ 0} \hfill & { - \theta_j \cdot \alpha _i } \hfill \\
 {\gamma_j \cdot \beta _i } \hfill & {\ \ 0} \hfill \\
\end{array} }} \right)|},\ (\forall j) \\
 \end{array}$$

$$\begin{array}{l}
0=d(\theta_l, \gamma_j)\geq \sum\limits_{i=1}^k {|\det
\left({{\begin{array}{*{20}c}
 {\theta_l \cdot \beta _j } \hfill & { - \theta_l \cdot \alpha _i } \hfill \\
 {\gamma_j \cdot \beta _i } \hfill & { - \gamma_j \cdot \alpha _i } \hfill \\
\end{array} }} \right)|}=\sum\limits_{i=1}^k {|\det
\left({{\begin{array}{*{20}c}
 {\ \ 0} \hfill & { - \theta_l \cdot \alpha _i } \hfill \\
 {\gamma_j \cdot \beta _i } \hfill & {\ \ 0 } \hfill \\
\end{array} }} \right)|},\ (\forall l\neq j) \\
 \end{array}$$

Since $\theta^T\cdot\alpha$, $\gamma^T\cdot\beta$ are both
nonsingular matrixes and $\theta_i\cdot\gamma_i=-1$, there is one
permutation $\sigma$ of $\{1,\cdots,k\}$ such that

$$\theta_i\cdot\alpha_{\sigma(i)}=\pm 1;\
\gamma_i\cdot\beta_{\sigma(i)}=\pm 1,\ (\forall i)$$

$$\theta_i\cdot\alpha_{\sigma(j)}=\gamma_i\cdot\beta_{\sigma(j)}=0,\ (\forall j\neq i)$$

Hence, we can obtain (18) and (19).\ {\bf QED}

\section{Right-hand Spheres and Left-hand Spheres in Simply Connected
3-Manifolds}

{\bf Theorem 4.}\ $(W\cup W^\prime;V_1,V_2), (W;V_1,V),
(W^\prime;V,V_2)$ satisfy the conditions of the Lemma 5. For every
$\beta_i$, we define a set $\beta_i(W,V)$ of the homotopy classes as
follows.
$$\beta_j(W,V)=\{[l]_V |l\subset V, l\sim \beta_i (in W)\}$$

Where $l\sim \beta_i\ (in W)$ denote that $l$ is homotopy
equivalence to $\beta_i$ in $W$. Let $G(\theta)$ generated by
$\{\theta_1,\cdots, \theta_k\}$ be a subgroup of $\pi_1(V)$. Then
$\{\theta_{\sigma(i)}^{\pm1}\}=\beta_i(W,V)\cap G(\theta)\
(i=1,\cdots,k)$, where $\sigma$ is a permutation of
$\{1,\cdots,k\}$.

{\bf Proof.}\ $W\cup W^\prime$ has a deformation retract

$$W\cup D_L^2(q_1)\cup \cdots\cup D_L^2(q_k)$$
where $\{D_L(q)\}$ are disjoint 2-discs, $W\cap D_L^2(q_i)=V\cap
D_L^2(q_i)=S_L^1(q_i)=\theta_i$.

Let $G(\theta)$ , generated by $\{\theta_1,\cdots, \theta_k\}$, be a
subgroup of $\pi_1(V)$, so each element $g\in G(\theta)$ is null
homotopy in $W\cup W^\prime$ and $W^\prime$.

$W^\prime$ has a deformation retract
$$V\cup D_L^2(q_1)\cup\cdots\cup D_L^2(q_k)$$

So we obtain the conclusions:

(1)\ In $W^\prime$, any closed path is homotopic onto $V$.

(2)\ In $V$, any closed path $l$ is null homotopy in $W^\prime$ if
and only if $[l]\in G(\theta)$.

If for some $i$, $\beta_i(W,V)\cap G(\theta)=\emptyset$, we will
show that $W\cup W^\prime$ is not simply connected.

Taking a positive number $\varepsilon$, $\varepsilon<1$, such that
two sets

$$X=W-f^{-1}(0),\ Y=f^{-1}(-\varepsilon,0]\cup D_L^2(q_1)\cup\cdots\cup
D_L^2(q_k)$$ are both path connected open subsets of $W\cup
D_L^2(q_1)\cup\cdots\cup D_L^2(q_k)$, $X\cap
Y=f^{-1}(-\varepsilon,0)$ is a open product manifold
$V\times(-\varepsilon, 0)$. Therefore, $\{X,Y,X\cap Y\}$ is a path
connected open covering of $W\cup D_L^2(q_1)\cup\cdots\cup
D_L^2(q_k)$. Since $\beta_i(W,V)\cap G(\theta)=\emptyset$, $\beta_i$
is not null homotopy in $Y$, moreover $\beta_i$ also is not null
homotopy in $X$. According to Van Kampen theorem, $W\cup
D_L^2(q_1)\cup\cdots\cup D_L^2(q_k)$, being the deformation retract
of $W\cup W^\prime$,  is not simply connected, thus $W\cup W^\prime$
also is not simply connected. However, $W\cup W^\prime$ is simply
connected, hence $\beta_i(W,V)\cap G(\theta)\neq\emptyset$.

According to (16) and (17), we have

$$\theta^T=-(\theta^T\cdot\alpha)\beta^T,\pmod{[\pi_1(T(k)),
\pi_1(T(k))]}$$
$$\gamma^T=(\gamma^T\cdot\beta)\alpha^T,\pmod{[\pi_1(T(k)),
\pi_1(T(k))]}$$

If $k=1$, $\pi_1(V)$ is a commutative group, so
$\theta=-(\theta\cdot\alpha)\beta$,
$\gamma=(\gamma\cdot\beta)\alpha$. Since $\theta\cdot\gamma=\pm1$
and $\alpha\cdot\beta=1$, we obtain
$(\theta\cdot\alpha)(\gamma\cdot\beta)=\pm1$, so $\theta=\pm\beta$,
$\gamma=\pm\alpha$. Hence, $d(S_L^1(q), S_R^1(p))=1$.

If $k\geq2$, $\pi_1(V)$ is not a commutative group.

According to Theorem 3, we have

$$\theta_i=\pm\beta_{\sigma(i)}\pmod{[\pi_1(V), \pi_1(V)]}$$

$$\gamma_i=\pm\alpha_{\sigma(i)}\pmod{[\pi_1(V), \pi_1(V)]}$$
where $\sigma$ is a permutation of $\{1,\cdots,k\}$

We can assume, by proper choice of orientations and serial numbers,
that $\sigma(j)=j\ (j=1,\cdots,k)$, and $-\theta^T\cdot\alpha=E_k$,
so we have

\begin{equation}
\theta^T=-(\theta^T\cdot\alpha)\beta^T=\beta^T,\pmod{ [\pi_1(T(k)),
\pi_1(T(k))]}
\end{equation}

Let $e\in\beta_i(W,V)\cap G(\theta)$, then $e\cdot\beta_j=0,\
(\forall j)$. $e$ can be expressed as

\begin{equation}
e=x_1^{\lambda_1}x_2^{\lambda_2}\cdots x_m^{\lambda_m},\ \
x\in\{\theta_1,\cdots, \theta_k\}
\end{equation}
where any two successive elements are different, and
$\lambda_1,\cdots, \lambda_m$ are nonzero integers.

As we know, $\pi_1(V)$ is the quotient of the free group on the
generators $\{\theta_1, \gamma_1,\cdots, \theta_k, \gamma_k\}$
modulo the normal subgroup generated by the element (see [2])
$$\prod_{i=1}^k[\theta_i, \gamma_i]$$
so $G(\theta)=\bigotimes_{i=1}^k G(\theta_i)$ is the free product of
$k$ infinite cycle groups $G(\theta_i)$, and therefore expression
(22) is unique.

$\theta_i$ can be expressed as

\begin{equation}
\theta_i=b_{i1}b_{i2}\cdots b_{in_i},\ \ b\in \{\alpha_1^{\pm1},
\beta_1^{\pm1}, \cdots, \alpha_k^{\pm1},\ \beta_k^{\pm1}\}
\end{equation}

The new product expression of $e$ is obtained from (22), (23)

\begin{equation}
e=(a_{11}a_{12}\cdots a_{1k_1})_1^{\lambda_1}(a_{21}a_{22}\cdots
a_{2k_2})_2^{\lambda_2}\cdots (a_{m1} a_{m2}\cdots
a_{mk_m})_m^{\lambda_m}
\end{equation}
$$a_{ij}\in \{\alpha_1^{\pm1},
\beta_1^{\pm1}, \cdots, \alpha_k^{\pm1},\ \beta_k^{\pm1}\}$$

$\pi_1(W)=\bigotimes_{i=1}^kG(\beta_i)$ is the free product of $k$
infinite cycle groups $G(\beta_i)$, and $\alpha_i=S_R^1(p_i)=BdD_R^2
(p_i)$, $(D_R^2(p_i)\subset W)$ is null homotopy in $W$,
$e\in\beta_i(W,V)$, so in $\pi_1(W)$, $e$ is uniquely expressed as

\begin{equation}
e=y_1^{\mu_1} y_2^{\mu_2}\cdots y_h^{\mu_h}\ \ y_i\in \{\beta_1,
\cdots, \beta_k\}
\end{equation}

where any two successive elements are different, $\mu_1, \cdots,
\mu_h$ are nonzero integers.

The expression (25) will be obtained by using following technique:

We delete all factors of $\{\alpha_1^{\pm1}, \cdots,
\alpha_k^{\pm1}\}$ from (24) and remain factors of
$\{\beta_1^{\pm1}, \cdots, \beta_k^{\pm1}\}$ with same order as they
in (24). Then (25) is obtained by simplification.

Since $e\in\beta_i(W,V)\cap G(\theta)$, hence $e$ is homotopy
equivalence to $\beta_i$ in $W$, so the expression (25) is just
$e=\beta_i^{\pm1}$.

On the boundary $V$, each $\theta_i$ can be expressed in the form
below

$$\theta_i |_V=A_i \beta_i B_i;  \ \  A_i,\ B_i\in[\pi_1(V),\pi_1 (V)]$$

In $W$, each $\theta_i$ can be expressed in the form below

$$\theta_i |_W=A_i|_W \beta_i B_i|_W; \ \ A_i|_W,\ B_i|_W\in[\pi_1 (W),\pi_1(W)]$$

We have the following equations

$$e|_V=(A_{i1} \beta_{i1} B_{i1})^{\lambda_{i1}} (A_{i2} \beta_{i2} B_{i2})^{\lambda_{i2}}\cdots
(A_{im_i}\beta_{im_i}B_{im_i})^{\lambda_{im_i}}$$

$$e|_W=(A_{i1}|_W\beta_{i1} B_{i1}|_W)^{\lambda_{i1}}(A_{i2}|_W\beta_{i2} B_{i2}|_W)^{\lambda_{i2}}\cdots
(A_{im_i}|_W\beta_{im_i} B_{im_i}|_W)^{\lambda_{im_i}}$$

Since $e\in\beta_i(W,V)\cap G(\theta)$£¬we have $e^{-1}\in\beta_i
(W,V)\cap G(\theta)$. In consideration of (21), we can assume
$\lambda_{i1}=1$ and

$$e|_W=(A_i|_W\beta_i B_i|_W)(A_{i2}|_W\beta_{i2}B_{i2}|_W)^{\lambda_{i2}}\cdots
(A_{im_i}|_W\beta_{im_i}B_{im_i}|_W)^{\lambda_{im_i}}$$

If $B_i|_W\neq 1$,then

$$(B_i|_W)^{-1}=(A_{i2}|_W\beta_{i2}B_{i2}|_W)^{\lambda_{i2}}\cdots
(A_{im_i}|_W\beta_{im_i}B_{im_i}|_W)^{\lambda_{im_i}}A_i$$

So $\theta_{i2}|_W,\cdots,\theta_{im_i}|_W$ are the proper factors
of $\theta_i|_W$, and $\theta_{i2}|_W\neq\theta_i|_W$. According to
(21),if $\theta_{i2}|_W$ is a proper factor of $\theta_i|_W$, then
$\theta_i|_W$ is not be the proper factor of $\theta_{i2}|_W$.
Similarly, if $A_i|_W\neq 1$£¬then $\theta_{im_i}|_W\neq\theta_i|_W$
and $\theta_{im_i}|_W$ is the proper factor of $\theta_i|_W$. if
$\theta_{i2}|_W$ is a proper factor of $\theta_i|_W$, we give the
relation $\theta_{i2}|_W<\theta_i|_W$. The relation $<$ define a
partial ordering of the set $\{\theta_1|_W,\cdots,\theta_k|_W\}$.
Suppose $\theta_i|_W$ be a minimal factor, then $A_i|_W=B_i|_W=1$,
and

$$e|_W=\beta_i(A_{i2}|_W\beta_{i2}B_{i2}|_W)^{\lambda_{i2}}\cdots
(A_{im_i}|_W\beta_{im_i}B_{im_i}|_W)^{\lambda_{im_i}}$$

According to (21), in the above formula, $\beta_h^{\pm 1}\ (h\neq
i)$ are all homogeneous. So if $m_i\geq2$, then the expression (25)
is not reduced to $\beta_i^{\pm1}$, hence, $m_i=1$. If $m_i=1$,
$\lambda_1\neq\pm1$, then the expression (25) is not reduced to
$\beta_i^{\pm1}$, hence $m_i=1$, $\lambda_1=\pm1$. So
$e=\theta_i^{\pm1}$ and $\{\theta_i^{\pm1}\}=\beta_i(W,V)\cap
G(\theta)$.\ {\bf QED}

{\bf Lemma 6.} Suppose that $(W; V_0,V)$ is a triad of the oriented
smooth compact 3-manifold, $V_0=T(k),\ (k\geq1)$ and $W$  has
exactly one critical point $q$  of type 2. Let $S$  be a
1-submanifold in $V_0$ and $S_L^1(q)\subset V_0$ be the left-hand
sphere of $q$.

If $d(S, S_L^1(q))>0$, in $V_0$, then the gradient image of any
closed path $s\subset V_0$ in the homotopy class $[S]$ is not a
closed path on $V$ , namely, $S$ is not homotopy onto $V$.

If $S_0\subset V$ is any closed path, then $S_0$ is homotopy
equivalence into $V_0$ in $W$. Let $S_0(W,V_0)$ denote all path
lifting from $S_0$ into $V_0$ in $W$, then for any closed path $g\in
S_0(W,V_0)$, $d(g, S_L^1(q))=0$ in $V_0$.

{\bf Proof.}\ We may assume that $S$  and $S_L^1(q)$ have exactly
$m=d(S, S_L^1(q))$ cross points $x_1, x_2, \cdots, x_m$. Take
disjoint curve segments $\{l_1, l_2, \cdots, l_m\}$ on $S$  such
that the curve segment $l_i$ pass the point $x_i$, so $l_i$ and
$S_L^1(q)$ intersect at one point $x_i$. Let $J(y)\ (y\in
V_0-S_L^1(q))$ denote the gradient curve via the point $y$ and
$J(y)\ (y\in S_L^1(q))$ denote the union of the gradient curve from
$y$ to $q$ and the right-hand disc $D_R^1(q)$. Let $Q$ be a subset
of $V_0$, define $J(Q)\cap V$ being the gradient image of $Q$ in
$V$. Let $l\subset V_0-S_L^1(q)$ be a continuous curve, then
$q\not\in J(l)$, so $J(l)\cap V$ is homeomorphic to $l$. As the set
$D_R^1(q)\cap V=S_R^0(q)$ has exactly two points, hence, $J(l_i)\cap
V$ is two disjoint curve segments and $J(S)\cap V$ is $m$ disjoint
curve segments. If $m=d(S, S_L^1(q))>0$, then, $J(S)\cap V$ is not a
closed path in $V$, moreover, $s_0\cap S_L^1(q)$\ $(s_0\in [S])$ has
at least $m$ points, so $J(s_0)\cap V$ is not a closed path in $V$.

Let $J(y)\ (y\in V-S_R^0(q))$ denote the gradient curve via the
point $y$ and $J(x)\ (x\in S_R^0(q))$ be the union of the gradient
curve from $q$ to $x$ and the left-hand disc $D_L^2(q)$. Let $E$ be
a subset of $V$, define $J(E)\cap V_0$ being the gradient image of
$E$ in $V_0$.

Taking a closed path $g\subset V-S_R^0(q)$, because $q\not\in J(g)$,
$J(g)\cap V_0$ is homeomorphic to $g$.

Let $S_0$ be a closed path in $V$. As $V$ is a connected 2-manifold
and the right-hand sphere $S_R^0(q)$ has exactly two points, so
there exists $S_1\in [S_0]$ satisfying $S_1\subset V-S_R^0(q)$,
$J(S_1)\cap V_0$ is homotopy equivalence to $S_1$ in $W$. Suppose
that $S_0\subset V_0$ is homotopy equivalence to $S\subset V$ in $W$
and $d(S_0,S_L^1(q))>0$ in $V_0$, according to the first conclusion,
it is impossible that $S_0$ is homotopy equivalence to $S$. Hence
$d(S_0,S_L^1(q))=0$ in $V_0$.\ {\bf QED}

Remark. The above results can be generalized to the case of more
than one critical point of index 2.

{\bf Theorem 5.}\ Suppose that $(W; V_1,V)$,\ $(W^\prime; V, V_2)$,\
$(W\cup W^\prime; V_1, V_2)$ are three oriented smooth compact
3-manifolds and that $W\cup W^\prime$,\ $V_1$,\ $V_2$ are all simply
connected. $f:(W\cup W^\prime; V_1, V_2)\rightarrow R^1$ is a Morse
function with the critical points $p_1, \cdots, p_k\ (k\geq1)$ of
type 1 and the critical points $q_1, \cdots, q_k$ of type 2, $f^{-1}
(-2)=V_1$,\ $f^{-1}(0)=V$,\ $f^{-1}(2)=V_2$,\ $f(p_i)=-1$,\
$f(q_i)=+1,\ (i=1,\cdots,k)$. If $\theta_i\cdot\alpha_i=\pm1\
(i=1,\cdots,k)$, then, there are at least a integer $i$, such that
$d(\theta_i, \alpha_h)=0\ (i\neq h)$.

{\bf Proof.}\ If $\theta_i\cdot\alpha_i=\pm1\ (i=1,\cdots,k)$,
according to Theorem 4, there are at least a integer $i$, such that

$$\{\theta_i^{\pm1}\}=\beta_i(W,V)\cap G(\theta)\ (i=1,\cdots,k)$$

Suppose that $i\neq h$ and $d(\theta_i, \alpha_h)>0$ on $V$.
According to Smale's conclusion ([1] p37-44), $(W\cup W^\prime; V_1,
V')$ can be expressed as

$$(W_1\cup W_2\cup W_3\cup W_4; V_1, V^\prime)$$
where $(W_1; V_1, V_2)$ has exactly one critical point $p_i$; $(W_2;
V_2, V)$ has exactly the critical points $\{p_j \mid j\neq i\}$; on
the same horizontal plane; $(W_3; V, V_4)$ has only one critical
point $q_i$; $(W_4; V_4, V_5)$ has exactly the critical points
$\{q_j \mid j\neq i\}$ on the same horizontal plane.

As $(W_1; V_1, V_2)$ has exactly one critical point $p_i$, so $V_2$
with genus 1 is an oriented 2-manifold and $\pi_1(W_1)$ with the
generators $\beta_i$ is the infinite cycle group, $\beta_i$ is a
1-submanifold in $V_2$.

Since $d(\theta_i, \alpha_h)>0$ on $V$, according to Lemma 6,
$\theta_i\ (\theta_i\subset V)$ is not homotopic onto $V_2$ in
$W_2$.

On the other hand, according to Lemma 6, $\beta_i$ can be homotopic
onto $V$ in $W_2$, all of the homotopy classes of $\beta_i$ on $V$
is just $\beta_i(W,V)$. If $[e]\in\beta_i(W,V)$, then $e$ is
homotopic onto $V_2$ in $W_2$. It is obtained from $\theta_i\cdot
\alpha_i=\pm1$, (21) and Theorem 4 that
$\{\theta_i^{\pm1}\}\subset\beta_i(W,V)$, so $\theta_i$ is homotopic
onto $V_2$ in $W_2$. Two contradictory conclusions show that
$d(\theta_i, \alpha_h)=0\ (i\neq h)$ is true.\ {\bf QED}

\section{Proof of Main Conclusion}

{\bf Lemma 7.}\ Assuming $T(k)$ is a differentiable, oriented and
closed 2-submanifold with the genus $k$; $M$ and $M^\prime$ are
smooth closed, transversely intersecting 1-submanifolds. Suppose
that the intersection numbers at $p, q\in M\cap M^\prime$ are +1 and
-1 respectively. Let $C$ and $C^\prime$ be the smoothly imbedding
arcs in $M$ and $M^\prime$ from $p$ to $q$. If $C$ and $C^\prime$
enclose a 2-disc $D$ (with two corners) with $Int D\cap (M\cap
M^\prime)=\emptyset$. Then, there exists an isotopy $h: T(k)\times
I\rightarrow T(k)$ such that

(1)\ \ $h_0$ is the identity map;

(2)\ \ The isotopy is the identity in a neighborhood of $M\cap
M^\prime-\{p, q\}$;

(3)\ \ $h_1(M)\cap M^\prime= M\cap M^\prime-\{p, q\}$.

{\bf Proof.}\ $M$ and $M^\prime$ are 1-manifolds, so there are two
one-sided collars $M\times[0,1)\subset T(k)$ with $M\times0=M$ and
$M^\prime\times[0,1)\subset T(k)$ with $M^\prime\times0=M^\prime$.
Take a small positive number $\varepsilon$ and two arcs $C_1$, $C_2$
with $C\subset C_1$,\ $C^\prime\subset C_2$. We may assume that
$C_1\times\varepsilon$ and $C_2\times\varepsilon$ transversely
intersect at two points $p^\prime$, $q^\prime$ and enclose a 2-disc
$E^\prime$ with $D\subset IntE^\prime$.

\begin{figure}[htb]
\begin{center}
\includegraphics[scale=.35]{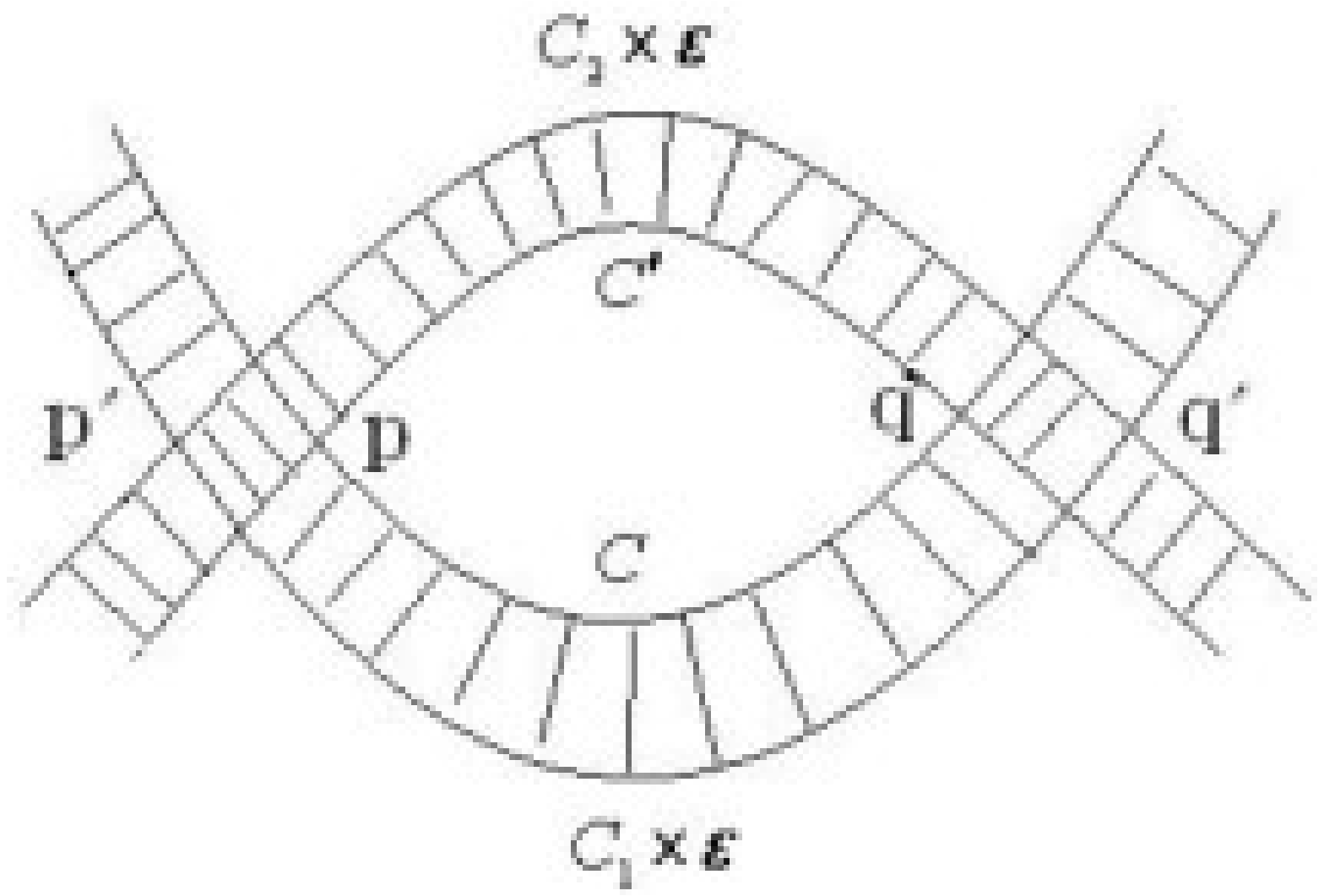}%
\end{center}
\end{figure}

Because the boundary
$BdE^\prime=(C_1\times\varepsilon)\cup(C_2\times\varepsilon)$ has
exactly two corners, by use of slight perturbation within a small
neighborhood of $p^\prime$, $q^\prime$, we obtain the smooth
boundary $S$ of a 2-disc $E$ with $D\subset IntE$.

Let $\Psi : E\rightarrow D^2$  denote a diffeomorphism of $E$ and
$D^2$.

Taking two arcs $L=M\cap E$  and $L^\prime=M^\prime\cap E$, then
$C\subset L$,\ $C^\prime\subset$,\ $L^\prime$; $L$ and $S=BdE$
transversely intersect at two points $\{x_1, x_2\}$; $L^\prime$ and
$S=BdE$ transversely intersect at two points $\{y_1, y_2\}$; $L$ and
$L^\prime$ only transversely intersect at two points $\{p, q\}$.

Since the intersection numbers of $M$ and $M^\prime$ at $p, q$  are
+1 and -1 respectively, there exists a line segment $\gamma$ in
$D^2$, $\gamma\cap S^1=\{a, b\}$, $a, b$ separate the boundary $S$
into two arcs $S_1$ and $S_2$ such that $\Psi(x_1), \Psi(x_2)\in
IntS_1$, $\Psi(y_1), \Psi(y_2)\in IntS_2$. Moreover, $\gamma$
separate $D^2$ into two closed area $A, B$ such that $\gamma\cup
S_1=Bd A$,\ $\gamma\cup S_2=Bd B$.

Taking a diffeomorphism  $\Omega:OD^2\rightarrow R^2$ with
$\Omega(\gamma)=R^1\times0$. We can assume, by proper choice of
$\Omega$, that $\Omega\circ\Psi(L)\cap R^-$ and
$\Omega\circ\Psi(L^\prime)\cap R^+$ are both the compact subset
(curve segments). Since $D\subset IntE$ is a compact subset,
$\Omega\circ\Psi(D)$ is a compact of $R^2$ and
$\Omega\circ\Psi(C\cup C^\prime)=Bd\Omega\circ\Psi(D)$.

Let $X_1(x)=(0, 1)$  be the unit vector field and
$\delta_r:R^2\rightarrow R^1$ be a differentiable function
satisfying the conditions
$$0\leq \delta_r(x)\leq1,\ \ \forall x\in R^2$$
$$\delta_r(x)=1,\ \ \forall x\in D_r^2$$
$$\delta_r(x)=0,\ \ \forall x\in R^2-D_{r+1}^2$$
where $D_r^2$ is a 2-disc with the radius $r$ in $R^2$.

$X_r(x)=\delta_r(x)X_1(x),\ \forall x\in R^2$ is a new vector field.
We claim $\Omega\circ\Psi(D)\subset OD_r^2$  by taking a large $r$.
The vector field determine one differentiable isotopy $h:R^2\times
R^1\rightarrow R^2$. Since $\Omega\circ\Psi(L)\cap R^-$ and
$\Omega\circ\Psi(L^\prime)\cap R^+$ are both the compact subset of
$R^2$, $h_{t_0}(\Omega\circ\Psi(L))\cap
\Omega\circ\Psi(L^\prime)=\emptyset$ can be obtained by taking
sufficient large $t_0$. According to the definition of $h_t$, we
obtain:
$$h_t(x)=x,\ \ \forall x\in R^2-D_{r+1}^2,\ \forall t\in R^+$$

The local isotopy $H_t=(\Omega\circ\Psi)^{-1}\circ
h_t\circ(\Omega\circ\Psi): IntE\times[0, t_0]\rightarrow IntE$ can
be extended to the total isotopy $H_t: T(k)\times[0, t_0]\rightarrow
T(k)$ satisfying the conditions
$$H_0\ \mbox{is the identity map}$$
$$H_t(x)=x,\ \ \forall x\in T(k)-E,\ \forall t\in[0, t_0]$$
$$H_{t_0}(M)\cap M^\prime=M\cap M^\prime-\{p, q\}$$

The desire isotopy is obtained.\ {\bf QED}

{\bf Lemma 8.}\ Let $\{\alpha_1, \alpha_2, \cdots, \alpha_k\}$ be
disjoint 1-submanifolds in $T(k)$ and $\theta$ be a 1-submanifold
satisfying the conditions $d(\theta, \alpha_j)=0\ (j=1,\cdots,h)$,
then, there exists finite diffeomorphisms which are isotopic to
identity, such that $\theta\cap\alpha_j=\emptyset\ (j=1,\cdots,k)$.

{\bf Proof.}\ If $\theta\cap \alpha_j\neq\emptyset$, we assume that
they are transversely intersect. Since $d(\theta, \alpha_j)=0$,
$\theta$ and $\alpha_j$ have cross points $\{x_{ij}, y_{ij}\mid i=1,
\cdots, r_j\}$, $\theta$ is separated into $2r_j$ smooth curve
segments $\{l_{j1}, \cdots, l_{j 2r_j}\}$ with $l_{jt}(0)\in
\{x_{ji}\mid i=1,\cdots,r_j\}$ and $l_{jt}(1)\in \{y_{ji}\mid
i=1,\cdots,r_j\}$, $\alpha_j$ is separated into $2r_j$ smooth curve
segments $\{g_{j1}, \cdots, g_{j 2r_j}\}$ with $g_{jt}(0)\in
\{x_{ji}\mid i=1,\cdots,r_j\}$ and $g_{jt}(1)\in \{y_{ji}\mid
i=1,\cdots,r_j\}$.

Since $d(\theta, \alpha_j)=0$, there exist $l_{ja}$ and $g_{jb}$
which enclose a 2-disc $D$ with
$IntD\cap(\theta\cap(\alpha_1\cup\cdots\cup\alpha_h))=\emptyset$.
According to lemma 7, there exists an isotopy $h_t:T(k)\times
I\rightarrow T(k)$, the isotopy is the identity in a neighborhood of
$\theta\cap (\alpha_1\cup\cdots\cup\alpha_h)-\{g_{jb}(0),
g_{jb}(1)\}$ and $h_1(\theta)\cap \alpha_j=\theta\cap
\alpha_j-\{g_{jb}(0), g_{jb}(1)\}$. Therefore, there exists finite
isotopies, and $\theta\cap \alpha_j=\emptyset\ (j=1,\cdots,)$.\ {\bf
QED}

{\bf Lemma 9.}\ Suppose that $(W; V_1, V)$,\ $(W^\prime; V, V_2)$,\
$(W\cup W^\prime; V_1, V_2)$ satisfy the conditions of theorem 5,
and $\theta_i\cdot\alpha_i=\pm1\ (i=1,\cdots,k)$, then $(W\cup
W^\prime; V_1, V_2)$ can be expressed as:
$$C_1 C_1^\prime C_2 C_2^\prime\cdots C_k C_k^\prime;\ \ p_i\in IntC_i,\ q_i\in
IntC_i^\prime$$ where $(C_i; V_i, U_i)$ and $(C_i^\prime; U_i,
V_{i+1})$ are the elementary cobordisms, and  $C_i C_i^\prime$,\
$V_j$ are all simply connected.

{\bf Proof.}\ According to theorem 4, we can assume
$\{\theta_1^{\pm1}\}=\beta_1(W,V)\cap G(\theta)$.

According to Theorem 5 and Lemma 7, Lemma 8, there are the isotopies
such that

\begin{equation}
S_L^1(q_1)\cap S_R^1(p_j)=\emptyset\ \ (j\geq 2)
\end{equation}

\begin{equation}
S_L^1(q_1)\cdot S_R^1(p_1)=\pm1
\end{equation}

Hence, we can alter the gradient field of $W^\prime$ satisfying (26)
and (27) in $V$.

On the basis of Smale's theorem on rearrangement of critical points
([1]), we have
$$W\cup W^\prime=C_1 C_1^\prime W_1 W_1^\prime;\  p_1\in C_1,\ q_1\in C_1^\prime,\
p_i\in W_1,\ q_i\in W_1^\prime\ (i\geq2)$$ And in $U_1$

$$S_L^1(q_1)\cdot S_R^1(p_1)=\pm1$$

Since $V_1$ is a simply connected 2-manifold and  $C_1$ has exactly
one critical point $p_1$ of type 1, that $U_1$ is a compact
2-manifold with the genus 1. $C_1^\prime$ has exactly one critical
 point $q_1$ of type 2, then the characteristic embedding associate
to $q_1$ is

$$\varphi: S^1\times OD^1\rightarrow U_1$$

So, $V_2=\chi(U_1, \varphi)$. $\chi(U_1, \varphi)$ denote the
quotient manifold obtained from the disjoint sum
$$(U_1-\varphi(S^1\times0))+(OD^2\times S^0)$$
by identifying

$$\varphi(u, \theta v)=(\theta u,v);\ \ u\in S^1,\ v\in S^0,\
0<\theta<1$$ where $\varphi(S^1\times0)$ is the left-hand sphere
$S_L^1(q_1)$ and $(0\times S^0)$ is the right-hand sphere
$S_R^0(q_1)$.

If $S_L^1(q_1)$ separate $U_1$ into two 2-manifolds with the
boundaries, then for any closed path $l$, the intersection number of
$l$ and $S_L^1(q_1)$ equal to zero, $l\cdot S_L^1(q_1)=0$, but
$S_L^1(q_1)\cdot S_R^1(p_1)=\pm1$, hence $S_L^1(q_1)$ dose not
separate $U_1$. So $V_2$ is connected.

In $(C_1^\prime; V_2, U_1)$, the critical point $q_1$ has index 1.
If the genus of $V_2$ is $g$, then the characteristic embedding
associate to $q_1$ is  $\varphi_R : S^0\times OD^2\rightarrow V_2$,
so the genus of $U_1=\chi(V_2, \varphi_R)$ is $g+1$. It is known
that the genus of $U_1$ is 1, hence $g=0$, $V_2$ is diffeomorphic to
$S^2$.

There is a path connected open covering $\{X, Y, X\cap Y\}$ of $C_1
C_1^\prime W_1 W_1^\prime$ satisfying $C_1 C_1^\prime\subset X$,
$W_1 W_1^\prime\subset Y$, $X\cap Y=S^2\times(-1,1)$. It is obvious
that $\pi_1(X\cap Y)=1$. Since $\pi_1(W\cap W^\prime)=1$,
$\pi_1(C_1C_1^\prime)=1$ and $\pi_1(W_2 W_2^\prime)=1$ can be
obtained by Van. Kampen theorem.

It is possible to alter the gradient vector field and rearrange the
critical points such that
$$W_1 W_1^\prime=C_2 C_2^\prime W_2 W_2^\prime;$$
$$\pi_1(C_2 C_2^\prime)=1,\ \ \pi_1(W_2 W_2^\prime)=1$$
$$V_2,\ V_3\ \mbox{are simply connected}$$

This procedure will continue until we derive the final conclusion.\
{\bf QED}

{\bf Lemma 10.}\ Suppose that $(W\cup W^\prime; V_1, V_2)$ is a
triad of the oriented compact 3-mnifold and $W\cup W^\prime$,\
$V_1$,\ $V_2$ are all simply connected. $(W; V_1, V)$ has exactly
one critical point $p$ of type 1, $(W^\prime; V, V_2)$ has exactly
one critical point $q$ of type 2. If in $V$, $S_L^1(q)\cdot
S_R^1(p)=\pm1$, then $(W\cup W^\prime; V_1, V_2)$ is a product
manifold $(S^2\times[0,1]; S^2\times0, S^2\times1)$.

{\bf Proof.}\ Since $V_1$ is diffeomorphic to $S^2$, $\lambda(p)=1$
and $\pi_1(W\cup W^\prime)=1$, $V$ is an oriented compact 2-manifold
and $V$ is diffeomorphic to $S^1\times S^1$. It is well-known that
$S^1\times S^1=R^1\times R^1/Z\times Z$ and $R^2$ is the universal
covering space of $S^1\times S^1$. Let $\Psi: R^2\rightarrow
S^1\times S^1$ be a covering mapping, $S_L^1(q)$, $S_R^1(p)$ are
both 1-submanifolds in $S^1\times S^1$. the intersection number of
$S_L^1(q)$ and $S_R^1(p)$ is +1 or -1, $S_L^1(q)\cdot
S_R^1(p)=\pm1$.

It is assumed that $S_R^1(p)$ has a path lifting  $L_0$ with the
origin as starting point,and the ending point of  $L_0$ is $(a, b)$,
where $a, b$ are two integer numbers. Then, all the path liftings of
$S_R^1(p)$ with $(na, nb)$ ($n$ are taken from all the integer
numbers) as the starting point to make up a smooth curve $L$ in
$R^2$. $L$ separates $R^2$ into two connected areas.

Since $S_L^1(q)\cdot S_R^1(p)=\pm1$, it is assumed that they are
transversely intersect with cross points $\{x_1, x_2,\cdots, x_r\}$
($r$ is an odd number).Taking one point $x\in\{x_1, x_2,\cdots,
x_r\}$ and one path lifting $L^\prime$ of $S_L^1(q)$ in $R^2$ such
that $L^\prime(0)\in \Psi^{-1}(x)\cap L$, then $L^\prime\cap L$ just
has $r$ points $\{x_1^\prime, x_2^\prime,\cdots, x_r^\prime\}\subset
R^2$. $\{x_1^\prime, x_2^\prime,\cdots, x_r^\prime\}$ separates
$L^\prime$ into $r-1$ curve segments $\{g_1, g_2,\cdots, g_{r-1}\}$
and separates $L$ into $r+1$ curve segments $\{l_1, l_2,\cdots,
l_{r+1}\}$.

If $r\geq3$, then there are two curve segments $l_\lambda$ and
$g_\mu$ which enclose a 2-disk such that $l_\lambda\cap g_\mu=\{x_j,
x_h\}$, and $(l_\lambda\cdot g_\mu)_{x_j}=1$, $(l_\lambda\cdot
g_\mu)_{x_h}=-1$. According to lemma 7, there exists an isotopy
$h_t,\ (0\leq t\leq1)$, the isotopy keeps the points nearer to
$L^\prime\cap L-\{x_j,x_h\}$ unmovable and $h(L^\prime)\cap
L=L^\prime\cap L-\{x_j,x_h\}$.

By use of the finite isotopies, $L^\prime$ and $L$ will have just
have one cross point. $\Psi: R^2\rightarrow S^1\times S^1=V$ may
bring these isotopies into $V$, enabling $S_L^1(q)$ and $S_R^1(p)$
to have only one cross point. According to Theorem 5.4 of [1] and
[3] (First cancellation theorem), it is possible to alter the
gradient vector field such that a new Morse function $f: W\cup
W^\prime \rightarrow R^1$ having no any critical point, so $(W\cup
W^\prime; V_1, V_2)$ is a product manifold $(S^2\times[0,1];
S^2\times0, S^2\times1)$.\ {\bf QED}

{\bf Lemma 11.}\ Let $(W\cup W^\prime; V_1, V_2)$ be an oriented
smooth 3-manifold; $W\cup W^\prime$,\ $V_1$,\ $V_2$ are all simply
connected. If there are $k$ critical points $\{p_1, p_2,\cdots,
p_k\}$, of type 1 in $(W; V_1, V)$ and they are on one same
horizontal plane. And if there are $k$ critical points $\{q_1,
q_2,\cdots, q_k\}$ of type 2 in $(W^\prime; V, V_2)$ and they are on
one same horizontal plane. Then, $(W\cup W^\prime; V_1, V_2)$ is
diffeomorphic to
$$(S^2\times[0,1]; S^2\times0, S^2\times1)$$

{\bf Proof.}\ According to Lemma 9 and Lemma 10, there exists a
Morse function $f: (W\cup W^\prime; V_1, V_2)\rightarrow R^1$
without any critical point, so $(W\cup W^\prime; V_1, V_2)$ is a
product manifold $(S^2\times[0,1]; S^2\times0, S^2\times1)$.\ {\bf
QED}

Let $M^3$ be a simply connected and compact smooth 3-manifold,
$BdM^3=\emptyset$.

According to the rearrangement theorem of cobordisms, we assume that
there is a Morse function of self-indexing such that

$$W_k=f^{-1}[k-\frac{1}{2}, k+\frac{1}{2}],\ \ (k=0,1,2,3)$$
$$V_{k+}=f^{-1}(k+\frac{1}{2})$$
$$f(p)=index(p),\ \mbox{at each critical point}\ p\ \mbox{of}\ f$$

{\bf Lemma 12.}\ $V_{1+}$ is the compact connected 2-manifold.

{\bf Proof.}  As $M^3$ is simple connection, so $V_{1+}$ is an
oriented 2-manifold. If $V_{1+}$ has $m$ connected components,
$V_{1+}=F_1+F_2+\cdots+F_m$; and every $F_i$ is an oriented, closed
2-manifold. All the critical points $p_1, p_2,\cdots, p_k$ of type 1
in $W_1$ are located on the same level $f^{-1}(1)$. All the critical
points $q_1, q_2,\cdots, q_n$ of type 2 in $W_2$ are located on the
same level $f^{-1}(2)$. So $W_1\cup W_2$ has a deformation retract
([1]Theorem 3.14)

$$D_R^2(p_1)\cup\cdots\cup D_R^2(p_k)\cup V_{1+}\cup D_L^2(q_1)\cup\cdots\cup
D_L^2(q_n)$$ where $D_R^2(p)$, $D_L^2(q)$ are 2-discs, and
$D_R^2(p)\cap V_{1+}=S_R^1(q)$, $D_L^2(q)\cap V_{1+}=S_L^1(q)$.

Therefore, $W_1\cup W_2$ just has $m$ connected components.

$W_0$ has exactly the critical points $o_1, o_2,\cdots,o_{n_1}$ of
type 0 locating on the same level $f^{-1}(0)$, so $W_0$ has a
deformation retract

$$D_R^3(o_1)\cup\cdots\cup D_R^3(o_{n_1})\cup V_{0+}$$

where $D_R^3(o_i)$ are disjoint 3-discs, and $D_R^3(o)\cap
V_{0+}=S_R^2(o)$.

$W_3$ has exactly the critical points $r_1, r_2,\cdots, r_{n_2}$ of
type 3 locating on the same level $f^{-1}(3)$, so $W_3$ has a
deformation retract
$$V_{2+}\cup D_L^3(r_1)\cup\cdots\cup D_L^3(r_{n_2})$$
where $D_L^3(r_j)$ are disjoint 3-discs, and $D_L^3(r)\cap
V_{2+}=S_L^2(r)$.

$W_0\cup W_1\cup W_2\cup W_3=M^3$ has a deformation retract
$$D_R^3(o_1)\cup\cdots\cup D_R^3(o_{n_1})\cup W_1\cup W_2\cup D_L^3(r_1)\cup\cdots
\cup D_L^3(r_{n_2})$$

Hence $M^3$ has exactly $m$ connected components. Since $M^3$ is
simple connection, $V_{1+}$ is the compact connected 2-manifold.\
{\bf QED}

{\bf Lemma 13.}\ Let $M^3$ be a simply connected and compact smooth
3-manifold, $BdM^3=\emptyset$. Then there exists a Morse function
$f: M^3\rightarrow R^1$, $f$ has exactly one critical point of type
0 and one critical point of type 3.

{\bf Proof.}\ Suppose that $f: M^3\rightarrow R^1$ has $n$ critical
points $o_1, o_2,\cdots, o_n$ of type 0 and $f$ is of self-indexing.
Then $W_0$ is just the sum of $n$ disjoint 3-discs
$W_0=D_1+D_2+\cdots+D_n$, $o_i\in D_i$ and
$V_{0+}=BdD_1+BdD_2+\cdots+BdD_n=S_R^2(o_1)+S_R^2(o_2)+\cdots+S_R^2(o_n)$.

Let $\varphi_i: S^0\times OD^2\rightarrow V_{0+}$ be the
characteristic embedding corresponding to the critical point $p_i$
of type 1, then the left-hand sphere
$\varphi_i(S^0\times0)=S_L^0(p_i)$ has exactly two points in
$V_{0+}$.

$W_1$ has a deformation retract
$$V_{0+}\cup D_L^1(p_1)\cup\cdots\cup D_L^1(p_k)$$
where $V_{0+}\cap D_L^1(p)=S_L^0(p)$ and $D_L^1(p_1),\cdots,
D_L^1(p_k)$ are disjoint 1-discs.

$V_{1+}$ is the compact connected 2-manifold, therefore $(W_1;
V_{0+}, V_{1+})$ is a connected 3-manifold.

If $n\geq2$, because $(W_1; V_{0+}, V_{1+})$ is connected, there is
$\varphi_j$ satisfying $\varphi_j(-1\times0)\in S_R^2(o_1)$ and
$\varphi_j(1\times0)\in S_R^2(o_2)$. Hence $S_R^2(o_1)$ and
$S_L^0(p_j)$ intersect at one point. According to the First
Cancellation Theorem, it is possible to alter the gradient vector
field $\xi$ of $f$ such that the gradient vector field $\xi^\prime$
of a new Morse function $f^\prime$, $f^\prime$  has exactly the
critical points $\{o_2,\cdots, o_n\}$ of type 0 and the critical
points $\{p_1,\cdots, \hat{p_j},\cdots, p_k\}$ of type 1, the other
critical points do not change.

This procedure will continue until we derive a Morse function with
one critical point of type 0. On the other hand, $-f:M^3\rightarrow
R^1$ has exactly the critical points $\{r_1,\cdots, r_{n_2}\}$ of
type 0, those critical points can be removed by using the same
technique till having one critical point of type 0. Hence there
exists a Morse function $f: M^3\rightarrow R^1$ which has exactly
one critical point of type 0 and one critical point of type 3.\ {\bf
QED}

{\bf Theorem 6.}\ A smooth compact simply connected 3-manifold $M^3$
is homeomorphic to $S^3$.

{\bf Proof.}\ We assume that $f: M^3\rightarrow R^1$ is the
self-indexing and $f$ has exactly one critical point of type 0 and
one critical point of type 3. Let $\chi(M^3)$ denote Euler
characteristic of $M^3$, it is well-known that $\chi(M^3)=0$.
According to Morse theorem ([4]), $f$ has exactly $k$ critical
points of type 1 and $k$ critical points of type 2. Moreover,
$W_1\cup W_2$,\ $V_{0+}$,\ $V_{2+}$ are all simply connected, so
$(W_1\cup W_2; V_{0+}, V_{2+})$ is a product manifold, there exists
a Morse function $f_0: M^3\rightarrow R^1$, $f_0$ has exactly one
critical point of type 0 and one critical point of type 3 and has
not any critical points of type 1, 2. Hence $M^3$  is homeomorphic
to $S^3$. Moreover, as every twist 3-sphere is diffeomorphic to
$S^3$ ([5], [6]) and $M^3$ is a twist 3-sphere, so $M^3$ is
diffeomorphic to $S^3$.

\end{document}